\input amstex
\documentstyle{amsppt}
\magnification=\magstep1
\NoBlackBoxes
\vsize=8.8truein
\loadmsbm
\def\cf{\Cal F}

\def\bs{\backslash}
\def\diam{\text{diam}}
\def\dist{\text{dist}}
\def\Vol{\text{vol}}
\NoRunningHeads
\TagsOnLeft

\topmatter
\title Foliation by Constant Mean Curvature Spheres on Asymptotically Flat
Manifolds
\endtitle
\author Rugang Ye
\endauthor
\affil Department of Mathematics\\ University of California\\
Santa Barbara, CA 93106\\  \&
 \\Mathematics Institute \\ Ruhr-University Bochum \\ 
\endaffil
\endtopmatter
\document

\subheading{\S 0 Introduction}

The main result of this paper is the following.

\proclaim{Main Theorem}  Let $M^{n+1}$ be an asymptotically flat manifold,
$n\geq
2$, and $\Omega$ an end of $M$ having nonzero mass.  Then there is on $\Omega$ a
smooth codimension one foliation $\cf_o$ by constant mean curvature spheres. 
$\cf_o$ is balanced and regular at $\infty$.  Moreover, it is the unique weakly
balanced and regular $C^2$ foliation on $\Omega$ by closed hypersurfaces of
constant mean curvature.
\endproclaim

``Balanced at $\infty$" means that near $\infty$ the leaves approach geodesic
spheres of a fixed center.  ``Weakly balanced" roughly means that the ``geodesic
centers" of the leaves do not shift to $\infty$ as fast as the farthest
points on the
leaves.  ``Regular at $\infty$" means that the rescaled second fundamental
form of
the leaves is uniformly bounded.  For the precise definitions we refer to
the next
section.  Note that $\cf_o$ actually foliates $\overline{\Omega\bs K}$ for a
compact region $K$.  In the statement of the theorem the phrase ``on
$\Omega$" is
used in a more general way than standard.  For the precise meaning of uniqueness
we refer to the Uniqueness Theorem in \S 2.

We obtained the existence part of this result (and a somewhat weaker uniqueness
result) at the end of 1988.  (In [Ye1],  we showed that around a
nondegenerate critical
point of the scalar curvature function in a Riemannian manifold, there exists a
unique regular foliation by constant mean curvature spheres.  We mentioned
that the
arguments extend to yield the existence result stated in the Main Theorem.)  Then
we obtained the uniqueness part of the Main Theorem.   Recently there has
been more
interest in this problem and several colleagues have inquired about the
details of
the proof of this result.  Also recently,  we showed in [Ye2] that in 
dimension $3$
all {\it diameter-pinched\/} \ (see Definition 5) $C^2$ foliations by constant mean curvature
spheres are
regular.  Since ``diameter-pinched" implies ``weakly balanced",  the Main
Theorem
implies the following

\proclaim{Strong Uniqueness Theorem}  Let $M$ be a 3-dimensional asymptotically
flat manifold, then on each end of nonzero mass of $M$ there is a unique
diameter-pinched $C^2$ foliation by constant mean curvature spheres. 
\endproclaim

For details we refer to [Ye2].  We remark that the ``diameter-pinched"
condition is
(much) weaker than the ``balanced" condition.

Asymptotically flat manifolds arise in general relativity.  The
significance of the
above results is to provide a canonical and regular geometric structure on
asymptotically flat manifolds.  Though the usual concept of asymptotical
flatness
is sufficient for applications, it is expressed in terms of coordinates and
hence is
not geometrically canonical.  (This is a problem on S. T. Yau's list [Ya]
of open
problems in differential geometry.  Note that in [B] a good understanding of the
asymptotical coordinates is provided.)  The canonical geometric structure
provided
by the above results serves to make the concept of asymptotical flatness more
geometrical and canonical.  In particular, it immediately makes the mass a
geometric invariant.  It also serves as a geometric linkage between different
asymptotical coordinates by the way of its construction.  Indeed, it should be
possible to completely characterize asymptotically flat ends of nonzero mass in
terms of balanced and regular foliations by constant mean curvature
spheres.  One
also expects further applications in general relativity.  A philosophical
implication
is a concept of ``center of universe" which may be defined as the region surrounded
by the foliations.

The program of constructing foliations by constant mean curvature spheres on
asymptotically flat manifolds was initiated by S. T. Yau. ( See [Ya] and [CY].)  He and G.
Huisken have
an independent proof of existence of foliations by constant mean curvature
spheres on
asymptotically flat ends of {\it positive\/} mass [HY1][HY2]. They apply the mean
curvature
flow to deform Euclidean spheres in asymptotical coordinates.  Positive
mass implies
a stability estimate which is employed to show convergence of the flow.  On
the other
hand, Yau and Huisken [HY2] showed that on a 3-dimensional asymptotically flat
manifold of
positive mass, the constructed foliation is the unique foliation by {\it
stable\/}
spheres of constant mean curvature. In comparison, our existence and uniqueness
results hold for both positive and negative mass and in all dimensions.  But we note that  our  
conditions on foliations for the uniqueness are rather different from 
that of Yau-Huisken. Their condition is stability as just mentioned. Our
conditions are  weak balance and regularity in general dimensions, and diameter pinching 
in dimension 3. 
These  have 
a different flavor than Yau
and Huisken's stability uniqueness. It is unknown whether uniqueness 
holds without any condition. In the local Riemannian situation, we do
have such a universal uniquness result in dimension 3, see [Y2]. But the 
situation of asymptotically flat ends is more subtle, see the discussion 
below.

Now we would like to discuss the proof of the Main Theorem.  There are several
delicate aspects here.  It is natural to attempt to perturb Euclidean spheres in
asymptotical coordinates in order to construct constant mean curvature
spheres.  But
a naive application of the implicit function theorem does not work, because the
linearized operator of the constant mean curvature equation in the problem has a
nontrivial kernel.  To resolve this difficulty, we apply the crucial idea
in [Ye1] of
moving centers.  Namely we first perturb the center of the asymptotical
coordinates
and then perform normal perturbation of the Euclidean spheres.  Both in
[Ye1] and here
the effect of the center perturbations is asymptotically degenerate.  We
remove this
degeneracy by carefully expanding the equation and balancing the center
perturbations against the normal perturbations.  On the other hand, it is
important
to control the magnitude of the center perturbations in order to retain the
foliation
property.  The situation of asymptotically flat ends is more delicate than
the local
picture in [Ye1], because the asymptotical flat structure deteriorates when the
center is shifted too far away.  This problem is even more serious for
uniqueness
than for existence, and it is the reason for the requirement of weak balance.  A
priori, a foliation by constant mean curvature spheres can differ
dramatically from
the one we constructed.  We have to obtain strong geometric control of the
leaves
in order to derive uniqueness.  Note that weak balance is a fairly weak
geometric
condition.  We think that it is necessary for uniqueness.

We acknowledge interesting discussions with G. Huisken.

\subhead \S 1 Moving Centers and Perturbation\endsubhead

Since we  deal with ends of asymptotically flat manifolds, it is convenient to
introduce the concept of asymptotically flat ends.  An asymptotically flat
manifold
is then a complete Riemannian manifold which is the union of a compact
region and
finitely many asymptotically flat ends.

\subhead Definition 1\endsubhead  Let $M$ be a complete Riemannian manifold of
dimension $n+1$ with $n\geq 2$.  Let $g$ be the metric of $M$.  A closed domain
$\Omega$ of $M$ is called an {\it asymptotically flat end} if there is a
coordinate
map from $\Omega$ to $\Bbb R^{n+1}\bs\overset o\to{\Bbb B}_{R_o}(o)$ for some
$R_o>0$ such that on this coordinate chart the metric $g$ can be written
$$
g_{ij}(x)=(1+\frac{\sigma}{r^{n-1}})\delta_{ij}+h_{ij}(x),\tag1.1
$$
where $r=|x|$, $\sigma$ is a constant and the $h_{ij}$'s satisfy 
$$
\align
&h_{ij}=O(\frac 1{r^n}), h_{ij,k}=O(\frac 1{r^{n+1}}),h_{ij,k\ell}=O(\frac
1{r^{n+2}}),\\
&h_{ij,k\ell m}=O(\frac 1{r^{n+3}}),
h_{ij,k\ell mm'}=O(\frac 1{r^{n+4}}),\endalign
$$
as $r\to\infty$ ($h_{ij,k}$ etc. denote partial derivatives, e.g.
$h_{ij,k}=\dfrac{\partial}{\partial x^k}h_{ij}$).  The constant $\sigma$ is
called
the {\it mass\/} or {\it energy\/} of $\Omega$.  (This differs from the usual
definition [LP] by a dimensional factor.)

Note that more general concepts of asymptotically flat manifolds have been
introduced in the literature, see e.g. [LP] and [Ye2].  But asymptotically
flat manifolds
as defined here are the most important.  (See e.g. [SY].  For technical
reasons we
required decay of up to the fourth order derivatives of $h_{ij}$. This 
same condition is also assumed in [HY].)  Hence we focus
on them in this paper.

\subhead Definition 2\endsubhead  Let $\cf$ be a foliation of codimension 1
on an
asymptotically flat end $\Omega$, whose leaves are all closed.  We say that
$\cf$
is {\it balanced at\/} $\infty$ or {\it balanced}, if there is a point
$p\in\Omega$
such that
$$
\frac{\dist(p,S)}{\diam(p,S)}\to 1 \text{ for } S\in\cf \text{ as }
\dist(p,S)\to\infty,
$$
where $\diam(p,S)=\underset{q\in S}\to\max\dist(p,q)$.  

\subhead Definition 3\endsubhead  Let $\cf$ be as above.  For $S\in\cf$, let
$\Omega_S$ be the open domain on the outside of $S$.  Put
$$
s(S)=\underset{p\in\Omega\bs\Omega_S}\to\max\frac{\dist(p,S)}{\diam(p,S)}.
$$
A point $p\in \Omega\bs\Omega_S$ is called a {\it geodesic center\/} of $S$, if
the maximum $s(S)$ is achieved at $p$.  We say that $\cf$ is {\it weakly
balanced}, if there are geodesic centers $p(S)$ of $S$ such that
$$
\limsup\frac{\dist(p_o,p(S))}{\diam(p(S),S)}<1 \text{ for } S\in\cf \text{ as
}\diam(p_o,S)\to\infty,
$$
where $p_o$ is a fixed point in $\Omega$.  We shall denote the above limit
by $b(S)$.

\subhead Definition 4\endsubhead  Let $\cf$ be as above.  Furthermore, assume
that  $\cf$ is of class $C^2$.  We say that $\cf$ is {\it regular at}
$\infty$ or {\it regular}, if
$$
\underset{\diam(S)\to\infty}\to\limsup \|A_S\|_{C^o(S)}\diam(S)<\infty\text{ for
}S\in\cf,
$$
where $A_S$ denotes the second fundamental form of $S$.

\subhead Definition 5\endsubhead  Let $\cf$ be  as above.
We say
that $\cf$ is {\it
diameter-pinched\/} at $\infty$ or diameter pinched, if
$$
\limsup\frac{diam(p_o,S)}{dist(p_o,S)}<\infty\text{ for } S
\in\cf\text{ as
}diam(p_o,S)\to\infty,
$$
where $p_o$ is a fixed point in $\Omega$.

Let $\Omega$ be an asymptotically flat end of dimension $n+1$, $n\geq 2$, whose
mass $\sigma$ is nonzero.  We can identify $\Omega$ with $\Bbb
R^{n+1}\bs\overset o\to{\Bbb B}_{R_o}$ for some $R_o>0$.   Without loss of
generality, we assume $R_o=1$.  To construct a balanced and regular foliation by
constant mean curvature spheres on $\Omega$, we apply the moving center method
in [Ye1].  Let $\nu$ denote the inward Euclidean unit normal of $S^n:\
=\partial\Bbb
B_1(o)$ and $\alpha_r$ the dilation $x\mapsto rx$ for $r>0$.  For $\varphi\in
C^2(S^n)$ and $\tau\in\Bbb R^{n+1}$ we define
$S^n_\varphi=\{x+\varphi(x)\nu(x):x\in S^n\}$ and
$S_{r,\tau,\varphi}=\alpha_r(\tau(S^n_\varphi))$, where the action of $\tau$ is
defined to be the translation by $\tau$.  Note that $S_{r,\tau,0}=\partial\Bbb
B_r(\tau)$ and $S^n_\varphi$ is an embedded $C^2$ surface if only
$\|\varphi\|_{C^1}\leq\varepsilon_o$ for some $\varepsilon_o\in(0,\dfrac 14)$. 
(We use the standard metric on $S^n$ unless otherwise stated.)  For $0<r<\dfrac
14$, $\|\varphi\|_{C^1}\leq\varepsilon_o$, $\tau\in\Bbb R^{n+1}$ and $x\in
S^n$ with
$|\tau|+r+\|\varphi\|_{C^o}\leq 1$, we put\newline
\centerline{ $H(r,\tau,\varphi)(x)=\dfrac 1r$ (the
inward mean curvature of the surface $S_{r,\tau,\varphi}$}
\centerline{ at $\dfrac
1r(\tau+x+\varphi(x)\nu(x)))$,} 
\noindent where of course we use the metric $g$ on $\Omega$.  
By (1.1), it is easy to see that $H(r,\tau,\varphi)$ extends to $r=0$ in a
$C^3$ fashion. 
We are going to compute $H(r,\tau,0)$.  Fix $r_o,\tau_o$ and $x_o\in S^n$. 
We need to
compute the mean curvature of $\partial\Bbb B_{1/r_o}(\tau_o/r_o)$ near
$y_o=(\tau_o+x_o)/ r_o$.  Choose coordinate vector fields $X_1,\cdots,X_n$ on
$\partial\Bbb B_{1/r_o}(\tau_o/r_o)$ near $y_o$ such that they are
orthonormal at
$y_o$ (with respect to the metric $g$).  We can extend them to
coordinate vector fields around $y_o$ in the following way:  $\widetilde
X_i((\tau_o+x)/r)=\dfrac{r_o}r X_i((\tau_o+x)/r_o)$.  We set
$Y((\tau_o+x)/r)=-x$. 
Then we have for $y=(\tau_o+x)/r_o$
$$
\align
\sum^n_{i=1}g(\nabla_{X_i}Y,X_i)\bigg|_y&=\frac 12\sum_{i=1}Yg(\widetilde
X_i,\widetilde X_i)\bigg|_y\tag1.2\\
&=nr_o+\frac 12\sum_{i=1}Yg(X_i,X_i)\bigg|_y,
\endalign
$$
where $\nabla$ denotes the Levi-Civita connection and $X_i$ is the euclidean
parallel translation of the original $X_i$ along the radial lines.  But 
$$
\sum^n_{i=1}g(X_i,X_i)\bigg|_y=g_{ii}(y)-g_{ij}(y)x^i x^j,
$$
where the summation convention is used on the right hand side for $1\leq
i$, $j\leq
n+1$.  Hence
$$
\sum^n_{i=1}g(\nabla_{X_i}Y,X_i)\bigg|_y=nr_o+\frac
12(g_{ii,k}(y)x^k-g_{ij,k}(y)x^i
x^j x^k).\tag1.3
$$
Next let $Z$ denote the inward unit normal of $\partial\Bbb B_{1/r_o}(\tau_o)$. 
$Z$ can be constructed explicitly from $Y$ by standard methods of linear
algebra. 
Indeed, there are linearly independent vectors $v_1(y),\cdots,v_n(y)$ which are
linear functions of $x$ and are orthogonal to $x$ in the Euclidean metric. 
Applying
the Gram-Schmidt procedure we obtain from the $v_i$'s orthonormal tangent
vectors $w_1(y),\cdots,w_n(y)$ at $y$.  Every $w_i$ can be written in the form
$F(g_{ij}(y),x)$ for a smooth function $F$ of $t_{ij}$, $1\leq i$, $j\leq
n+1$ and
$x$.  Finally, set $Z_1=\sum\limits^n_{i=1}h(Y,w_i)w_i$, where $h$ denotes the
tensor $h_{ij}$ at $y$.  Then 
$$
Z=\frac{Y-Z_1}{g(Y-Z_1,Y-Z_1)^{1/2}}.
$$
We compute
$$
g(Y-Z_1,Y-Z_1)=1+\frac\sigma{|y|^{n-1}}+h_{ij}Y^i Y^j-2g_{ij}Z^i_1
Y^j+g_{ij}Z^i_1
Z^j_1,
$$
whence
$$
\align
Z-Y=-\frac\sigma{2|y|^{n-1}} Y&+\frac 1{|y|^{2n-2}}F_1
(g_{ij},h_{ij},|y|^{-1},x)\tag1.4 \\
&+h_{ij}F_{ij}(g_{k\ell},h_{k\ell},|y|^{-1},x)
\endalign
$$
with smooth functions $F_1$ and $F_{ij}$.  Next we observe that the covariant
derivative $\nabla$ can be explicitly written and that we can replace $X_i$ by
$w_i$ in (1.3).  Since the mean curvature of $\partial\Bbb B_{1/r_o}(\tau_o)$ is
given by $\sum\limits^n_{i=1}g(\nabla_{w_i}Z,w_i)$, we deduce from (1.3) that
$$
\align
&H(r,\tau,0)(x)=n(1-\frac\sigma{2|y|^{n-1}})+\frac 1{2r}(g_{ii,k}x^k-g_{ij,k}
x^ix^jx^k)\cdot\tag1.5\\
&(1-\frac\sigma{2|y|^{n-1}})
+\frac 1{r|y|^{n+1}}\widetilde F_1(g_{ij},
|y|g_{ij,k},h_{ij},|y|h_{ij,k},|y|^{-1},|y|^{-1}y,x)\\
&+\frac 1{|y|^n}\widetilde F_2(g_{ij},h_{ij},|y|^{-1},x)+\frac 1r
h_{ij,k}F_{ijk}(g_{\ell m},h_{\ell m},|y|^{-1},x)\\
&+\frac 1{r|y|} h_{ij}\widetilde F_{ij}(g_{k\ell},|y| g_{k\ell},h_{k\ell},|y|
h_{k\ell},|y|^{-1},|y|^{-1}y,x)\\
&+ h_{ij}\overline F_{ij}(g_{k\ell},h_{k\ell},|y|^{-1},x).
\endalign
$$
Here $r_o,\tau_o$ have been replaced by $r$ and $\tau$, $y=(\tau+x)/r$, the
variable for $g_{ij},g_{ij,k}$ etc. is $y$, and $\widetilde F_1,\widetilde F_2$,
$F_{ijk}, \widetilde F_{ij},\overline F_{ij}$ are smooth functions.  We have
$$
g_{ij,k}=\frac{\sigma(1-n)}{|y|^{n+1}} y^k\delta_{ij}+h_{ij,k},\tag1.6
$$
$$
|y|^{-1}=r|\tau+x|^{-1}=r(1-x^i\tau^i+\tau^i\tau^j b_{ij}(\tau,x))\tag1.7
$$
for smooth functions $b_{ij}$ which are defined for all $x\in S^n$ and
$\tau\in\Bbb
R^{n+1}$ with $|\tau|<1$.  From (1.5), (1.6) and (1.7) we then deduce $$
\align
H(r,\tau,0)(x)=n-\frac{\sigma n^2}2 r^{n-1}&+\frac{\sigma n(n-1)(n+1)}2
r^{n-1}\tau^i x^i\tag1.8\\
&+r^{n-1}\tau^i\tau^j\tilde b_{ij}(\tau,x)\\
&+r^{n-1}f(r,\tau,x),
\endalign
$$
with smooth functions $\tilde b_{ij}$, $f$ for $r>0$, $|\tau|<1$ and $x\in
S^n$ such
that 
\medskip
\noindent 1) $\|f(r,\tau,\cdot)\|_{C^1(S^n)}\leq C  (|\tau|) r,$\newline
\medskip
\noindent 2) $\|d_\tau f(r,\tau,\cdot)\|_{C^1(S^n)}\leq C (|\tau|) r,$\newline
\medskip
\noindent 3) $\|\dfrac{\partial f}{\partial r}(r,\tau,\cdot)\|_{C^1(S^n)}\leq
C(|\tau|)$\newline
\medskip

\noindent for a positive continuous function $C(t)$ defined for $|t|<1$,
where $d_\tau
f$ denotes the differential of $f$ in $\tau$.  We put $\tilde f=f+\tilde
b_{ij}\tau^i \tau^j$.

Our goal is to find solutions $\tau,\varphi$ of the equation 
$$
H(r,\tau,\varphi)=n-\frac{\sigma n^2}2 r^{n-1}.\tag1.9
$$
We consider $H(r,\tau,\cdot)$ as a mapping from $C^{2,1/2}(S^n)$ into
$C^{0,1/2}(S^n)$ and let $H_\varphi$ denote the differential of $H$ w.r.t.
$\varphi$. 
Set $g_{r,\tau}=r^2g$.  Then $H(r,\tau,\varphi)$ is the mean curvature of
$S^n$ and
$H_\varphi(r,\tau,0)$ is just the Jacobi operator $\Delta
+\|A\|^2+Rc(\overline\nu)$
on $S^n$ relative to the metric $g_{r,\tau}$, where $Rc$ denotes Ricci curvature
and $\overline\nu$ the inward unit normal of $S^n$.  We indicate the
dependence on
$g_{r,\tau}$ as follows 
$$
H_\varphi(r,\tau,0)=\Delta_{r,\tau}+\|A_{r,\tau}\|^2+Rc_{r,\tau}.\tag1.10
$$
Note that $g_{\tau,r}$ converges to the Euclidean metric in the $C^4$
topology as
$r\to 0$.  It follows that $H_\varphi(0,\tau,0)=L:=\Delta_{S^n}+n$, where
$\Delta_{S^n}$ is the standard Laplace operator on $S^n$.  Let $Ker$ denote the
kernel of $L$ which is spanned by the Euclidean coordinate functions.  We
have the
orthogonal $L_2$-decompositions $C^{2,1/2}(S^n)=Ker\oplus Ker^\bot$ and
$C^{0,1/2}(S^n)=Ker\oplus L(Ker^\bot)$.  Let $P$ denote the orthogonal
projection
from $C^{0,1/2}(S^n)$ onto $Ker$ and $T:Ker\to\Bbb R^{n+1}$ the isomorphism
sending $x^i\big|_{S^n}$ to $\bold e_i=$ the $i$-th coordinate basis.  Put
$\widetilde P=TP$.  Then 
$$
\widetilde P(H(r,\tau,0))=\frac 12
n(n-1)(n+1)\sigma\omega_{n+1}r^{n-1}\tau+r^{n-1}\widetilde P(\tilde
f(r,\tau,\cdot)),\tag1.11
$$
where $\omega_{n+1}=\Vol(\Bbb B_1(o))$.  To solve (1.9) we introduce the
following
expansions
$$
\align
H(r,\tau,\varphi)&=H(r,\tau,0)+r\int^1_0 \int^1_0 H_{\varphi
r}(sr,\tau,t\varphi)\varphi ds\ dt\tag1.12\\
&+\int^1_0\int^1_0 tH_{\varphi\varphi}(0,\tau,st\varphi)\varphi\varphi ds\ dt\\
&+L\varphi,
\endalign
$$
where the subscript $r$ means the partial derivative in $r$ and 
$$
H_{\varphi\varphi}(r,\tau,\psi)\varphi\varphi'=\frac
d{dt}H_{\varphi}(r,\tau,\psi+t\varphi')\varphi\big|_{t=0}\ .
$$

We first consider the equation  
$$
\widetilde P(H(r,\tau,r^{n-1}\varphi))=0\ .\tag1.13
$$
Dividing it by $r^{n-1}$ and applying (1.11), (1.12) we can reduce it to
the following
$$
\frac 12 n(n-1)(n+1)\sigma\omega_{n+1}\tau+\widetilde P(\tilde
f(r,\tau,\cdot))+r\widetilde P(q(r,\tau,\varphi))=0,\tag1.14
$$
where
$$
\align
q(r,\tau,\varphi)&=\int^1_0\int^1_0 H_{\varphi
r}(sr,\tau,tr^{n-1}\varphi)\varphi ds\
dt\\
&+r^{n-1}\int^1_0\int^1_0 tH_{\varphi\varphi}(0,\tau,st
r^{n-1}\varphi)\varphi\varphi ds\ dt.
\endalign
$$
(Note that $\widetilde P L=0$.)  Easy computations show that $q$ has the
following
properties\newline
\medskip
\noindent 1)
$\|q(r,\tau,\varphi)\|_{C^{0,1/2}(S^n)}\leq\beta(\|\varphi\|)$,\newline
\medskip
\noindent 2) $\|d_\tau
q(r,\tau,\varphi)\|_{C^0(S^n)}\leq\beta(\|\varphi\|)$,\newline
\medskip
\noindent 3) $\|q_\varphi(r,\tau,\varphi)\|\leq\beta(\|\varphi\|)$ with
$q_\varphi:C^{2,1/2}(S^n)\to C^{0,1/2}$,\newline
\medskip
\noindent 4) $\|\dfrac{\partial q}{\partial
r}(r,\tau,\varphi)\|_{C^{0,1/2}(S^n)}\leq\beta(\|\varphi\|)$.
\medskip
\noindent Here and in the sequel, $\beta$ denotes a positive continuous
function on
$\Bbb R$ and $\|\varphi\|=\mathbreak\|\varphi\|_{C^{2,1/2}(S^n)}$.

Let $Q(r,\tau,\varphi)$ denote the left hand side of (1.14).  Then for each
given bound
on $\varphi$, $Q(r,0,\varphi)$ approaches zero uniformly as $r\to 0$.  This
follows
from the properties of $q$ and $\tilde f$.  Also by virtue of these
properties, the
differential $d_\tau Q$ is nearly the identity for small $r$ and $\tau$. 
Hence we can
apply the implicit function theorem to obtain a unique solution
$\tau(r,\varphi)$ of
the equation (1.14) and hence (1.13) for small $r$ which lies near zero. 
There holds
$\tau(r,\varphi)\to 0$ as $r\to 0$.  We also have the following estimates 
$$
\|\tau_\varphi\|\leq\beta(\|\varphi\|)r,\quad
\fracwithdelims||{\partial\tau}{\partial
r}\leq\beta(\|\varphi\|).
$$
These estimates are easy consequences of the properties of $q$ and $\tilde
f$.  It
follows that
$$
|\tau(r,\varphi)|\leq\beta(\|\varphi\|)r.
$$

Now we replace $\tau$ in (1.9) by $\tau(r,\varphi),\ \varphi$ by
$r^{n-1}\varphi$, and
divide (1.9) by $r^{n-1}$.  Then (1.9) is reduced to 
$$
L\varphi+\frac{\sigma n(n-1)(n+1)}2 \tau^i(r,\varphi) x^i+\tilde
f(r,\tau(r,\varphi),\cdot)+rq(r,\tau(r,\varphi),\varphi)=0.\tag1.15
$$
By the above argument for finding $\tau(r,\varphi)$ and the properties of
$\tilde f,q$
we can find a unique solution $\varphi(r)$ of (1.15) for small $r$ which
lies near
zero.  There holds $\varphi(r)\to 0$ as $r\to 0$.  Moreover, we have the
following
estimates:
$$
\|\varphi(r)\|\leq Cr,\bigg|\bigg|\frac{d\varphi}{dr}\bigg|\bigg|\leq C
$$
for a constant $C>0$.

The uniqueness property of $\tau(r,\varphi)$ and $\varphi(r)$ can be stated
precisely
as follows

\proclaim{Proposition 1}  For each positive number $C$ there are positive
numbers
$R_1=R_1(C)$ and $R_2=R_2(C)$ with the following properties:  1) if
$Q(r,\tau,\varphi)=0$,\ $0<r\leq R_1$,\ $|\tau|\leq R_2$ and
$\|\varphi\|\leq C$, then
$\tau=\tau(r,\varphi)$; 2) if
$H(r,\tau(r,\varphi),r^n\varphi)=n-\dfrac{\sigma n^2}2
r^{n-1}$, $0<r\leq R_1$ and $\|\varphi\|\leq C$, then $\varphi=\varphi(r)$.
\endproclaim

We omit the easy proof.  (For part 2) one utilizes the equation (1.15) to
show that
$\varphi$ is small.)

Now we consider the family of surfaces $\cf_0=\{\widetilde
S_r=S_{1/r,\tau(r,\varphi(r)),r^{n-1}\varphi(r)}:0<r\leq r_o\}$, where
$r_o$ is chosen
as follows.  Set $C=\underset{r\leq R_1(1)}\to\sup\|\varphi(r)\|+1$,
$r_1=\max\{r:|\tau(r,\varphi)|\leq R_2(C)$ with $\|\varphi\|\leq C\}$ and then
$r_o=\min\{r_1,R_1(1),R_1(C)\}$.  These surfaces are diffeomorphic to $S^n$.  By
construction, $\widetilde S_r$ has constant mean curvature
$nr-\dfrac{\sigma n^2}2
r^n$.  Since all the equations and known functions in the above construction are
smooth, we conclude that $\cf_o$ is a smooth family of constant mean curvature
spheres.  Geometrically, we obtained this family by moving the center (the
origin) of
the Euclidean spheres $\partial\Bbb B_r(o)$ to $\tau(r,\varphi(r))/r$ and then
performing the Euclidean normal perturbation $r^{n-2}\varphi(r)$.  It remains to
show that $\cf_o$ is a foliation.  Put
$\psi(r,x)=r^{-1}\tau(r,\varphi(r))+r^{-1}x+r^{n-2}\varphi(r)(x)$, where
$x\in S^n$. 
By the estimates for $\tau$ and $\varphi$ we easily see that the maps
$v(r,\cdot):=\psi(r,\cdot)/|\psi(r,\cdot)|$ approach the identity map from
$S^n$ to
$S^n$ in the $C^2$ norm as $r\to 0$.  Hence $v(r,\cdot)$ is a smooth
diffeomorphism for small $r$.  We set $\tilde\varphi(r,x)=|\psi(r,v^{-1}(r,x))|$
where $v^{-1}(r,\cdot)$ denotes the inverse of $v(r,\cdot)$.  Then $\widetilde
S_r$ is the Euclidean normal graph of $\tilde\varphi(r,\cdot)$ over $S^n$. 
Employing
the estimates for $\tau$ and $\varphi$ we compute 
$$
\align
\frac{\partial\psi}{\partial r}&=-\frac x{r^2}+\frac
1r\frac{\partial\tau}{\partial
r}+\frac
1r\tau_\varphi\fracwithdelims(){d\varphi}{dr}-\frac\tau{r^2}+(n-2)r^{n-3}
\varphi(r)+r^{n-1}\frac{d\varphi}{dr}\\
&=-\frac 1{r^2}(x+O(r)),\\
&d_x\psi=\frac 1r(I+O(r)),\ d_x\,v=I+O(r),
\endalign
$$
and finally
$$
\frac{\partial\widetilde\varphi}{\partial r}=\frac
1{|\psi|}\psi\cdot\left(\frac{\partial\psi}{\partial
r}+(d_x\psi)\fracwithdelims(){\partial v^{-1}}{\partial r}\right)=-\frac
1{r^2}(1+O(r)). $$
We conclude that $\widetilde\varphi(r,x)$ is strictly decreasing w.r.t. $r$
for $r$
small (or better, $\widetilde\varphi(\rho^{-1},x)$ is strictly increasing
for $\rho$
large).  Hence $\widetilde S_r,\widetilde S_r$, are disjoint for small
$r,r'$ with
$r\neq r'$.  We replace $r_o$ by a smaller positive number if necessary. 
Then $\cf_o$
constitutes a smooth foliation.  It is easy to see that $\cf_o$ is balanced
and regular
at $\infty$.  Thus we have proven the existence part of the Main Theorem.

\subhead\S 2 Uniqueness\endsubhead

Let $\Omega$ be the asymptotically flat end considered above and $\cf_o$ the
constructed foliation.  Let $\Omega_1$ be the support of $\cf_o$, i.e. the
subdomain
of $\Omega$ foliated by $\cf_o$.  Note that $\Omega\bs\Omega_1$ is compact.  

\subhead Definition 3\endsubhead  A {\it constant mean curvature
foliation\/} is a
codimension one $C^2$ foliation with closed leaves of constant mean curvature.  

We have

\proclaim{Uniqueness Theorem}  Let $\cf$ be a weakly balanced and regular
constant
mean curvature foliation on a closed subdomain $\Omega'$ of $\Omega$ such that
$\Omega\bs\Omega'$ is compact.  Then either $\cf_o$ is a restriction of $\cf$ or
$\cf$ is a restriction of $\cf_o$.  In other words, $\cf_o$ is the unique
maximal
weakly balanced and regular constant mean curvature foliation in
$\Omega_1$.
\endproclaim

\demo{Proof} Let $\cf$ be as described in the theorem.  It is easy to see
that the
leaves of $\cf$ can be parametrized as a $C^2$ family $S_t$, $0<t\leq 1$ with
$S_t\neq S_{t'}$ if $t\neq t'$ and $\underset{t\to 0}\to\lim\diam
S_t=\infty$ (see
Lemma 2.1 in [Ye1]).  Then $S_t$ lies on the interior side of $S_{t'}$,
provided that
$t>t'$.  We put $\ell(t)=(\diam S_t)^{-1}$ and
$S^*_t=\alpha_{\ell(t)}(S_t)$, i.e.
$S^*_t$ is the dilation of $S_t$ by the factor $\ell(t)$.  Let
$g_t=\ell(t)^2\alpha^*_{1/\ell(t)}g$.  Then $g_t$ converges to the
Euclidean metric in
the $C^4$ topology away from the origin.  Each $S^*_t$ has constant mean
curvature in
the metric $g_t$.  Because $\cf$ is regular, the second fundamental form of
$S^*_t$
in $g_t$ is bounded above by a constant independent of $t$.

We claim that $\dist(o,S^*_{t_k})$ is uniformly bounded away from zero. 
Assume the
contrary.  Then we can find a sequence $S^*_{t_k}$ with $t_k\to 0$ which
converges
to an immersed surface $S_\infty$ in $\Bbb R^{n+1}\bs\{o\}$ of constant mean
curvature.  $S_\infty$ has uniformly bounded second fundamental form and
diameter
$1$.  The only possible self-intersections of $S_\infty$ are of the type that an
embedded piece of $S_\infty$ meets another from one side.  A priori,
$S_\infty$ may
have many components.  Let $S'_\infty$ be one component.  Then $\overline
S\,'_\infty\bs S'_\infty=\{o\}$.  Since $S_\infty$ has uniformly bounded second
fundamental form and constant mean curvature, the origin is a removable
singularity,
i.e. $\overline S\,'_\infty$ is a smooth immersed surface of constant mean
curvature. 
Applying the classical Alexandrov reflection principle we deduce that
 $\overline S\,'_\infty$ is a round sphere.  It follows that $\overline
S_\infty$ is
either a round sphere containing the origin or the union of two round
spheres meeting
tangentially at the origin.  But the condition of weak balance does not
allow such
limits as one readily sees.  Thus we arrive at a contradiction and the
claim is proven.

Consider an arbitrary sequence $S^*_{t_k}$ with $t_k\to 0$.  By the above
claim and
the above arguments, a subsequence converges to a round sphere $S_\infty$ of
diameter $1$ such that the origin lies in the open ball bounded by
$S_\infty$.  We
conclude that
$$
\widetilde H(t)\diam S_t=2n+h(t)\tag2.1
$$
for a function $h(t)$ which converges to zero as $t\to 0$, where
$\widetilde H(t)$
denotes the mean curvature of $S_t$ (in the metric $g$).  In particular
there is for
small $t$ a unique solution $r(t)$ of the equation
$$
\widetilde H(t)=nr-\frac{\sigma n^2}2 r^n.\tag2.2
$$
Now we set $\widehat S_t=\alpha_{r(t)}(\widetilde S_t)$ and $\widehat
g_t=r(t)^2\alpha^*_{1/r(t)}g$.  Then the mean curvature $H(t)$ of $\widehat
S_t$ in
$\widehat g_t$ is $n-\dfrac{\sigma n^2}2 r(t)^{n-1}$.  By the above arguments
concerning $S^*_t$ we deduce that every sequence $\widehat S_{t_k}$ with $t_k\to
0$ contains a subsequence converging to a round sphere $S_\infty$ such that the
origin lies in the open ball bounded by $S_\infty$.  Because of (2.1), $\diam
S_\infty=2$.  Let $a_\infty$ denote the center of $S_\infty$.  Then 
$$
|a_\infty|\leq b(\cf)<1.
$$
We conclude that for small $t$ every $\widehat S_t$ is the normal graph of
a smooth
function $\varphi_t$ over a round sphere of radius $1$ and center $a(t)$
such that
$|a(t)|\leq b(\cf)$ and $\varphi_t$ converges to zero in the $C^4$ norm as
$t\to 0$.  By
the arguments in the proof of Lemma 2.3 in [Ye1] we can find (for small $t$)
$v(t)\in\Bbb R^{n+1}$ such that $\underset{t\to 0}\to\limsup |v(t)|\leq b(\cf)$,
$\widehat S_t=S_{1,v(t),\widehat\varphi(t)}$ for a smooth function
$\widehat\varphi(t)$ over $S^n$ with $\widehat\varphi(t)\to 0$ in the $C^4$
norm as
$t\to 0$.  Moreover, the projection $P(\widehat\varphi(t))$ vanishes.  Roughly
speaking, we can move the center $a(t)$ slightly to achieve that the
defining function
of $\widehat S_t$ as a graph has zero kernel component.  Note that
$S_t=S_{1/r(t),v(t),\widehat\varphi(t)}$ and 
$$
H(r(t),v(t),\widehat\varphi(t))=n-\frac{\sigma n^2}2 r(t)^{n-1}.\tag2.3
$$

Applying (2.3), (1.8), (1.12) and the inequality $\underset{t\to
0}\to\limsup|v(t)|\leq b(\cf)<1$ we obtain for small $t$
$$
\|\widehat\varphi(t)\|\leq Cr(t)^{n-1}\tag2.4
$$
for a constant $C>0$.  Next we estimate $v(t)$.  To this end we apply (1.5) and
calculate $H(r,\tau,0)$ in a way somewhat different from (1.8).  We have
$$
\align
H(r,\tau,0)(x)&=n-\frac{n\sigma}2 r^{n-1}\frac 1{|x+\tau|^{n-1}}-
\frac{n(n-1)\sigma}2
r^{n-1}\frac 1{|x+\tau|^{n+1}}\tag2.5\\
&-\frac{n(n-1)\sigma}2 r^{n-1}\frac{\tau^k x^k}{|x+\tau|^{n+1}}+ r^{n-1}
f(r,\tau,x),
\endalign
$$
where $f$ is the same function as appearing in (1.8).  It follows from
(2.3), (2.4), (2.5)
and (1.12) that
$$
\bigg| P\left(\frac 1{|x+v(t)|^{n-1}} +(n-1)\frac
1{|x+v(t)|^{n+1}}+(n-1)\frac{v(t)\cdot x}{|x+v(t)|^{n+1}}\right)\bigg|\leq
Cr(t)\tag2.6
$$
for a positive constant $C$.  For a fixed small $t$ we change coordinates
such that
$v(t)=\ell\bold e_1$ for some $\ell\in(0,1)$.  Then
$$
\align
&-\int\limits_{S^n}\frac{x^1}{|x+v(t)|^{n-1}}d\Vol-(n-1)\int\limits
_{S^n}\frac{x^1}{|x+v(t)|^{n+1}}
d\Vol\\
&-(n-1)\int\limits_{S^n}\frac{(v(t)\cdot x)x^1}{|x+v(t)|^{n+1}} d\Vol\\
&=-\int\limits_{S^n}\frac{x^1}{|1+\ell x^1|^{n-1}}
d\Vol-(n-1)\int\limits_{S^n}\frac{x^1}{|1+\ell x^1|^n}d\Vol\cdot\\
&\geq(n-1)\int\limits_{x^1>0} x^1\left(\frac 1{|1-\ell x^1|^n}-\frac 1{|1+\ell
x^1|^n}\right) d\Vol\\
&\geq 2n(n-1)\int\limits_{x^1>0}\frac{(x^1)^2}{(1+\ell x^1)^n(1-\ell
x^1)^n}d\Vol\\
&\geq (n-1)(n+1)\omega_{n+1}\ell,
\endalign
$$
whence
$$
|v(t)|\leq Cr(t)\tag2.7
$$
for a constant $C>0$.

Setting $\varphi^*(t)=r(t)^{1-n}\widehat\varphi(t)$ we have
$Q(r(t),v(t),\varphi^*(t))=0$.  Applying (2.4), (2.7) and Proposition 1 we
deduce that
$v(t)=\tau(r(t),\varphi^*(t))$ for small $t$.  Applying (2.3) and
Proposition 1 we
then conclude that $\varphi^*(t)=\varphi(r(t))$ for small $t$.  Consequently
$S_t=\widetilde S_{r(t)}$ for small $t$.  A continuity argument then shows that
either $\cf_o$ is a restriction of $\cf$ or $\cf$ is a restriction of
$\cf_o$.\quad\qed
\enddemo

 \end